\documentclass[preprint,12pt]{elsarticle}
\makeatletter
\def\ps@pprintTitle{%
  \let\@oddhead\@empty
  \let\@evenhead\@empty
  \let\@oddfoot\@empty
  \let\@evenfoot\@oddfoot
}
\makeatother

\usepackage{amsmath}
\usepackage{amssymb}
\usepackage{amsthm}
\newtheorem{theorem}{Theorem}[section]
\newtheorem{lemma}[theorem]{Lemma}

\usepackage{enumitem}
\usepackage{upgreek}
\usepackage{mathrsfs}

\newtheorem{example}[theorem]{Example}
\usepackage{hyperref}
\usepackage[thinc]{esdiff}
\newtheorem{remark}{Remark}

\journal{}

\begin{document}

\begin{frontmatter}%



\title{Drazin and group invertibility in algebras spanned by  two idempotents}


\author[inst1]{Rounak Biswas}

\affiliation[inst1]{organization={Department of Mathematical and Computational Science, National Institute of Technology Karnataka},
            city={Surathkal}, 
            postcode={575025},
            state={Karnataka},
            country={India}}

\author[inst1]{Falguni Roy}

\begin{abstract}
For two given idempotents $p\text{ and }q$ from an associative algebra $\mathcal{A},$ in this paper, we offer a comprehensive classification of algebras spanned by the idempotents $p\text{ and }q$. This classification is based on the condition that $p\text{ and }q$ are not tightly coupled and satisfies $(pq)^{m-1}=(pq)^{m}$ but $(pq)^{m-2}p\neq (pq)^{m-1}p$ for some $m(\geq2)\in\mathbb{N}.$ Subsequently, we categorized all the group invertible elements and established an upper bound for Drazin index of any elements in these algebras spanned by $p,q$. Moreover, we formulate a new representation for the Drazin inverse of $(\alpha p+q)$ under two different assumptions, $(pq)^{m-1}=(pq)^m$ and $\lambda(pq)^{m-1}=(pq)^m,$ here $\alpha$ is a non-zero and $\lambda$ is a non-unit real or complex number. 

\end{abstract}

\begin{keyword}
Finite dimensional algebra \sep group inversion \sep Drazin inversion \sep idempotent

\MSC 15A30 \sep 15A09 \sep 16S15 \sep 17C27  
\end{keyword}

\end{frontmatter}
Email: rounak.207ma005@nitk.edu.in; royfalguni@nitk.edu.in
\section{Introduction}\label{Preliminary }
Throughout this paper, $\mathcal{A}$ always represent an associative algebra over the field of scalars $\mathbb{C}$ or $\mathbb{R},$ with an identity $I$. An element $a\in \mathcal{A}$ is said to be idempotent if $a^2=a$.  In this paper $p$ and $q$ will consistently refer to two idempotent elements in $\mathcal{A}$. For given $p,q$; alg$(p,q)$ denotes the subalgebra of $\mathcal{A}$ spanned by $p\text{ and }q$, i.e. alg$(p,q)$ is composed of all possible finite linear combinations of elements from the list \begin{align}\label{list}
    &p,pq,pqp,(pq)^2\cdots(pq)^m\cdots\\\notag
    &q,qp,qpq,(qp)^2\cdots(qp)^m\cdots.
\end{align}
Here, the order of an element from (\ref{list}) is defined by the number of factors present. The notations $(pq)_i,\text{ }(qp)_i$ denotes the upper and lower element of the $i$-th column in the list (\ref{list}).
For given $p,q$ if two members from the list (\ref{list}) with order difference $\leq1$ are equal, then following \cite{bottcher2013classification}, $p,q$ termed as tightly coupled idempotents. In simpler terms, we say that $p$ and $q$ are tightly coupled if, for some $i\in\mathbb{N}$, either the product $(pq)_i$ or $(qp)_i$ matches one of the five neighbouring elements from the list $(\ref{list})$. Specifically, this means that $(pq)_i$ (or $(qp)_i$) coincides with any of the five adjacent elements( called neighbourhood): ${(pq)_{i-1}, (pq)_{i+1}, (qp)_{i-1}, (qp)_i, (qp)_{i+1}}$ (or ${(qp)_{i-1}, (qp)_{i+1}, (pq)_{i-1}, (pq)_{i}, (pq)_{i+1}}$), for some $i\in\mathbb{N}$.  B{\"o}ttcher and Spitkovsky \cite{bottcher2013classification,bottcher2011certain} studied the algebra spanned by $p$ and $q,$ when they are tightly coupled.
Depending on the relation between $p$ and $q,$ 
they classified alg$(p,q)$ into four different types $Z_m,U_m,D_m$, $D_m^*$ and established \begin{align}\label{idk}
    D_m^*&\cong Z_{m-2}\oplus D_2^*\text{ } (m\geq 2),\notag \\
    D_m&\cong Z_{m-2}\oplus D_2\text{ } (m\geq 2),\\
    U_m&\cong Z_{m-1}\oplus U_1 \text{ } (m\geq 1),\notag
    \end{align}where $m\in\mathbb{N}.$
Here $Z_m=\mathrm{alg }(p,q),$ when $p,q$ satisfy one of the following conditions depending on $m$\begin{align*}
    \begin{cases}
        (qp)_k=(pq)_k=0,
\text{ if }m=2k-2\\
(qp)_k=0\text{ or }(pq)_k=0,\text{if }m=2k-1
\end{cases},
\end{align*} and for construction of the algebras $U_m,D_m^*,D_m$ interested readers can refer to \cite{bottcher2013classification}.
However, when attempting to extend this classification beyond tightly coupled scenarios, B{\"o}ttcher and Spitkovsky \cite{bottcher2013classification,bottcher2011certain} acknowledged the challenges and complexities involved in categorizing the algebra alg$(p,q)$. In this article, we establish that starting from one of the two elements in any column of the list (\ref{list}), there is only one way to go beyond tightly coupled.  Specifically, if one initiates with $(pq)^m$(or $(qp)^m$) from the $2m$-th column of the list (\ref{list}), then the only way go beyond tightly coupled is when $(pq)^m=(pq)^{m-k}$(or$(qp)^m=(qp)^{m-k}$), for some $k\in\mathbb{N}$ and $k<m.$ Furthermore, we provide a complete classification of the algebra alg$(p,q)$, when $p\text{ and }q$ are not tightly coupled and satisfies $(pq)^{m-1}=(pq)^{m},$ $(pq)^{m-2}p\neq (pq)^{m-1}p.$
 The investigation of these algebraic structures is motivated by the goal of comprehending Drazin and group inevitability within this context.

For $a\in \mathcal{A},$ the Drazin inverse of $a$ is the unique element $b\in \mathcal{A}$(denoted by $a^D$), satisfying
\begin{equation}\label{drazin_def}
    ab=ba,\text{ }ab^2=b\text{ and }a^{k+1}b=a^k,
\end{equation}for some $k\in\mathbb{N}\cup \{0\}$. The least $k$ satisfying (\ref{drazin_def}), is known as the Drazin index $i(a)$ of $a$. In the particular case $i(a)=1,$ the Drazin inverse is termed as group inverse( symbolized by $a^g$). As demonstrated by Drazin \cite{drazin1958pseudo}, the Drazin invertibility of $a$ is equivalent to existence of two elements $x,y\in \mathcal{A}$ satisfying \begin{equation}\label{Drazin_Def_2}
    a^{k_1+1}x=a^{k_1}\text{ and }ya^{k_2+1}=a^{k_2},
\end{equation} for some positive integer $k_1,k_2$. The smallest $k_1,k_2$ satisfying (\ref{Drazin_Def_2}) is referred to as the left and right index of $a$, respectively. Moreover if $k_1,k_2$ is finite and the least integer satisfying (\ref{Drazin_Def_2}) then $k_1=k_2=i(a)$ with $a^D=a^{k_1}x^{k_1+1}=y^{k_1+1}a^{k_1}.$ One of the interesting problems in Drazin inverse theory is to investigate Drazin invertibility of $\alpha p+\beta q,$ where $\alpha, \beta$ are scalars. This specific problem has garnered great interest from various researchers in recent years. Deng \cite{deng2009drazin} considered this problem for idempotent operators on a Hilbert space, when $\alpha,\beta\in\{1,-1\}$ which is extended in the setting of a Banach algebra by Zhang, Wu \cite{zhang2012drazin}, where they provided an expression of $(\alpha p+\beta q)^D$ assuming $pqp=0$ and $pqp=pq$, for non zero scalars $\alpha, \beta$. More results on this problem can be found in \cite{shi2013drazin,chen2014drazin,xie2012drazin,chen2020drazin}. In particular, for tightly coupled idempotents, $p,q$; alg$(p,q)$ becomes finite-dimensional, and hence, every element in this algebra becomes Drazin invertible. However, the same statement is not true in the case of group inverse. 
B{\"o}ttcher and Spitkovsky \cite{bottcher2012group} solved this problem of group invertibility in $Z_m$, and that is enough to settle the case for other algebras due to (\ref{idk}). But we observe that the Theorem 2.2. \cite{bottcher2012group}, which is  essential for investigating group invertibility in $Z_m$, lacks consideration of the case when $Z_m$ and $\mathcal{A}$ possess different identity elements. An example (Example \ref{example}) narrating this observation will be discussed in Section \ref{g&d}. This paper addresses this gap by presenting a comprehensive version of Theorem 2.2. \cite{bottcher2012group}, and utilizing this, we solve the problem of group invertibility in alg$(p,q),$ when $p,q$ are not tightly coupled and satisfies $(pq)^{m-1}=(pq)^{m},$ $(pq)^{m-2}p\neq (pq)^{m-1}p$ for some $m\in\mathbb{N}.$ Providing a representation of the Drazin inverse (or group inverse) for any element in associative algebra poses a significant challenge once it is established that the element possesses such an inverse. In this paper, we will provide a new representation of $(\alpha p+q)^D$ under two different conditions; the first one is when $(pq)^{m-1}=(pq)^{m}$ and another one is when $\lambda (pq)^{m-1}=(pq)^m,$ here $\alpha$ is a non-zero and $\lambda$ is a non-unit scalar.

The paper is organized in the following manner.  Section \ref{finite}, presents the classification of algebras generated by $p\text{ and }q$, when $(pq)^m=(pq)^{m-1}$. Section \ref{g&d}, deals with the Darzin index and the group invertible elements in the algebras discussed in section \ref{finite}. Representations of the Drazin inverse of $(\alpha p+q)$ under two different assumption $(pq)^{m-1}=(pq)^m$ and $\lambda (pq)^{m-1}=(pq)^m$ are in Section \ref{representation}.
 \section{Algebras spanned by  non-tightly coupled idempotents}\label{finite}
Due to the complexity in extending the neighbourhood beyond tightly coupled,  B{\"o}ttcher and Spitkovsky  \cite{bottcher2013classification} settled the case for the algebra alg$(p,q)$, when $pqp=p$. But it turns out that there is only one way to extend the neighborship beyond tightly coupled; our next lemma establishes this.
\begin{lemma}
    For $m,k\in\mathbb{N}$, where $k< m$, if $p$ and $q$ satisfy one of the following conditions, then they become tightly coupled:
    \begin{enumerate}[label=(\roman*)]
        \item $(pq)^{m-k}p=(pq)^m;$
        \item $(qp)^{m-k}q=(pq)^m;$
        \item $(qp)^{m-k}=(pq)^m.$
    \end{enumerate}
    \begin{proof}
        \begin{enumerate}[label=(\roman*)]
            \item Since $(pq)^{m-k}p=(pq)^m$, multiplying both sides by $q$ from the right yields $(pq)^{m-k+1}=(pq)^m$. Consequently,  we get $(pq)^{m-k}p=(pq)^{m-k+1}$, indicating that $p$ and $q$ are tightly coupled.
            \item Multiplying $(pq)^m=(qp)^{m-k}q$ by $p$ from left we obtain $(pq)^m=(pq)^{m-k+1}$, therefore $(qp)^{m-k}q=(pq)^{m-k+1}.$
            \item Similarly multiplying $(pq)^m=(qp)^{m-k}$ by $p$ from left we obtain $(qp)^{m-k}=p(qp)^{m-k}.$
        \end{enumerate}
    \end{proof}
\end{lemma}
Therefore, starting from $(pq)^m$ for some $m\in\mathbb{N}$, there is only one way to extend the neighbourship beyond tightly coupled, i.e. when $(pq)^m=(pq)^{m-k}$ for some integer $k<m.$ For this paper, our focus is on the particular case when $k=1$.  In the following lemma, we provide some properties of the idempotents $p,q$ when they satisfy $(pq)^m=(pq)^{m-1}$ for some $m\in\mathbb{N}$.
\begin{lemma}\label{n=n-k}
    Let $p,q$ satisfies $(pq)^{m-1}=(pq)^m$ for $n\in\mathbb{N}$,
    \begin{enumerate}[label=(\roman*)]
        \item $(pq)^{m-1}=(pq)^{m+k}$, for $k\in\mathbb{N}$;
        \item if $(pq)^{n_1}=(pq)^{n_2}$, where $n_1,n_2\in\mathbb{N}$ such as $1\leq n_2<n_1\leq m$ then $(pq)^m=(pq)^{m-i}$ for some $i\in\mathbb{N}$;
        \item if $(pq)^m=(pq)^{m-k}$ for some $k\in\mathbb{N}$, where $1<k<m$, then ${pq}^{m-2}=(pq)^{m-1}.$
    \end{enumerate}
    \begin{proof}
        \begin{enumerate}[label=(\roman*)]
            \item One can verify this easily.
            \item If $n_1=m$ then $i=m-n_2$. Now if $n_1<m$ then multiplying both side of $(pq)^{n_1}=(pq)^{n_2}$ with $(pq)^{m-n_1}$ we get $$(pq)^m=(pq)^{n_2+m-n_1}.$$ Hence $(pq)^m=(pq)^{m-i}$ where $i=n_1-n_2.$
            \item Since $(pq)^m=(pq)^{m-k}$ then multiplying both side by $(pq)^{k-2}$ we obtain 
            \begin{align*}
                (pq)^{m-2}=(pq)^{m+k-2}{}=(pq)^{m-1}.
            \end{align*}
        \end{enumerate}
    \end{proof}
\end{lemma}
\begin{remark}
    According to Lemma \ref{n=n-k}, if $m$ is the least positive integer satisfying $(pq)^m=(pq)^{m-1}$ then $(pq)^j\neq(pq)^l$ for $j,l\in\mathbb{N}$ where $j<l< m$.
\end{remark} 
In the remaining part of this section, we will provide a complete classification of all algebras spanned by the idempotents $p,q;$ where $p,q$ are not tightly coupled and satisfies $(pq)^{m-1}=(pq)^{m}$ but $(pq)^{m-2}p\neq (pq)^{m-1}p$ for some $m\in\mathbb{N}.$ Because of these assumptions on $p,q$, the infinite list $(\ref{list})$ terminate to a finite list; hence, alg$(p,q)$ becomes finite dimensional.  Before proceeding further, let's revisit the following lemma from \cite{bottcher2013classification}; it provides the complete classification of  alg$(p,q)$ when $pqp=p$. Note that in the context of any algebra $\mathcal{X}$ and for $A\in\mathcal{X}$; the notations $\mathcal{N}(\mathcal{X})$ and $\sigma(A)$, represent the collection of nilpotent elements in $\mathcal{X}$ and the spectrum of $A$ within $\mathcal{X}$, respectively.
\begin{lemma}\cite{bottcher2013classification}
Upto isomorphism, there are precisely four algebras, namely $W_3,W_3\oplus Z_1,W_4,$ $W_4\oplus Z_1,$ spanned by $p,q$;  where $pqp=p$ and $q,p,qp,pq$ are pairwise distinct.
    \begin{enumerate}[label=(\roman*)]
        \item For the case when $q=qpq$, $q+p=qp+pq,$ we have $\mathrm{alg}(p,q)\cong W_3$ where $\mathrm{dim}\text{ } W_3=3$, $\mathrm{dim}\text{ }\mathcal{N}(W_3)=2.$
        \item If $ q\neq qpq$ with $qp+pq=qpq+p$, we have $\mathrm{alg}(p,q)\cong  W_3\oplus Z_1$ where  $\mathrm{dim}\text{ }(W_3\oplus Z_1)=4$ and $\mathrm{dim}\text{ }\mathcal{N}(W_3\oplus Z_1)=2.$
        \item For the case when $q=qpq$ but $q+p\neq qp+pq$, we have $\mathrm{alg}(p,q)\cong  W_4,$ with $\mathrm{dim}\text{ }W_4=4$,  $\mathrm{dim}\text{ }\mathcal{N}(W_4)=3.$
        \item If $ q\neq qpq$ with $ qp+pq\neq qpq+
        p$, then $\mathrm{alg}(p,q)\cong  W_4\oplus Z_1$ where $\mathrm{dim}\text{ }(W_4\oplus Z_1)=5$, $\mathrm{dim}\text{ }\mathcal{N}(W_4\oplus Z_1)=3$.
    \end{enumerate}
\end{lemma}
\begin{lemma}\label{inde}
    Let $p\text{, }q$ are not tightly coupled and satisfy 
    \begin{equation}\label{least}
       (pq)^{m-1}=(pq)^{m} \text{ but }(pq)^{m-2}p\neq (pq)^{m-1}p,
    \end{equation}
     for some $m(\geq 2)\in\mathbb{N},$ where $m$ is the least positive integer satisfying $(\ref{least})$. Then the set
     \begin{equation}\label{collection}
         \{p,q,pq,qp,\cdots,(pq)^{m-2}p,(qp)^{m-2}q\}
     \end{equation} is linearly independent.
     \begin{proof}
         Since $m$ is the least positive integer satisfying (\ref{least}), therefore by Lemma \ref{n=n-k}, $(pq)^{m-1}\neq (pq)^{m-k}$ for $k\in\mathbb{N}$ where $2\leq k<m.$ Now consider an element from the linear span of the set $(\ref{collection})$ which is claimed to be zero,
         \begin{equation}\label{linear_com}
             \sum_ix_i(pq)_i+\sum_iy_i(qp)_i,
         \end{equation} where $x_i,y_i$ are scalars. Let $l$ denote the smallest order of the products that appear in $(\ref{linear_com})$ with a nonzero coefficient. First let assume $l$ to be an even integer, then if $x_l\neq 0$, multiplying (\ref{linear_com}) by $p$ in left side and with $q$ in right side, respectively, we obtain,
         \begin{equation}\label{linear_com_2}
            x_l(pq)^{\frac{l}{2}}+\sum_{i\in\{l+1,l+3,\cdots\}}(x_i+x_{i+1})(pq)^{\frac{i+1}{2}}+\sum_{i\in\{l,l+2,\cdots\}}(y_i+y_{i+1})(pq)^{\frac{i+2}{2}}=0. 
         \end{equation}
      Next multiplying (\ref{linear_com_2}) by $(pq)^{m-\frac{l}{2}-2}$ we get\begin{equation}\label{linear_com3}
         x_l(pq)^{m-2}+\left(\sum_{i>l}x_i+\sum_{i\geq l}y_i\right)(pq)^{m-1}=0.
     \end{equation} Now again multiplying (\ref{linear_com3}) by $pq$ and subtracting from (\ref{linear_com3}) we get $$x_l\left((pq)^{m-2}-(pq)^{m-1}\right)=0,$$hence we obtain $x_l=0$ as $(pq)^{m-2}-(pq)^{m-1}\neq 0,$ a contradiction. Similarly, let $y_l\neq0,$ then multiply (\ref{linear_com}) by $pq$ from left and from right by $p(qp)^{m-\frac{l}{2}-2}$, hence we have \begin{equation}\label{linear_com4}
         y_lp(qp)^{m-2}-\sum_{i\neq l}(x_i+y_i)p(qp)^{m-1}=0.
     \end{equation} Since $p(qp)^{m-2}\neq p(qp)^{m-1},$ from (\ref{linear_com4}) we obtain $y_l=0$, a contradiction.
     
     In the case when $l$ is odd, in a similar way, one can verify that $x_l=y_l=0$. Hence, the collection $(\ref{collection})$ consists of linearly independent elements.
     \end{proof}
\end{lemma}
If the assumptions outlined in Lemma \ref{inde} are satisfied by $p$ and $q$, then in an analogous way like \cite[Lemma 5.1]{bottcher2013classification}, using Lemma \ref{inde} one can verify the subsequent lemma.
\begin{lemma}\label{alg_class}
    Let $p\text{ and }q$ satisfy the assumptions of Lemma $\ref{inde}$.
    \begin{enumerate}[label=(\roman*)]
        \item If $(qp)^{m-1}= (qp)^{m}$ and $(qp)^{m-1}+(pq)^{m-1}= (qp)^{m-1}q+(pq)^{m-1}p$ then $$p,q,pq,qp,\cdots ,(pq)^{m-1},(qp)^{m-1}q$$ forms a linearly independent set where $\mathrm{dim}\text{ }\mathrm{alg}(p,q)=4m-3.$
        \item If $(qp)^{m-1}\neq (qp)^{m}$ and $(qp)^m+(pq)^{m-1}= (qp)^{m-1}q+(pq)^{m-1}p$ then $$p,q,pq,qp,\cdots, (pq)^{m-1}p,(qp)^{m-1}q$$ forms a linearly independent set where $\mathrm{dim}\text{ }\mathrm{alg}(p,q)=4m-2.$
        \item If $(qp)^{m-1}= (qp)^{m}$ and $(qp)^{m-1}+(pq)^{m-1}\neq (qp)^{m-1}q+(pq)^{m-1}p$ then $$p,q,pq,qp,\cdots, (pq)^{m-1}p,(qp)^{m-1}q$$ forms a linearly independent set where $\mathrm{dim}\text{ }\mathrm{alg}(p,q)=4m-2.$
        \item If $(qp)^{m-1}\neq (qp)^{m}$ and $(qp)^m+(pq)^{m-1}\neq (qp)^{m-1}q+ (pq)^{m-1}p$ then $$p,q,pq,qp,\cdots, (pq)^{m-1}p,(qp)^{m-1}q,(qp)^m$$ forms a linearly independent set with $\mathrm{dim}\text{ }\mathrm{alg}(p,q)=4m-1.$

    \end{enumerate}
\end{lemma}
Now, we are ready to establish the main result of this section.
\begin{theorem}
    Let $p,q$ satisfy the conditions of Lemma $\ref{inde}$.
    \begin{enumerate}[label=(\roman*)]
        \item \label{kk}If $(qp)^{m-1}= (qp)^{m}$ and $(qp)^{m-1}+(pq)^{m-1}= (qp)^{m-1}q+(pq)^{m-1}p$ then here $\mathrm{alg}(p,q)\cong Z_{4m-6}\oplus W_3$ with $\mathrm{dim}\text{ }\mathcal{N}\left(\mathrm{alg}(p,q)\right)=4m-6.$
        \item \label{3k}In the case when $(qp)^{m-1}\neq (qp)^{m}$ and $(qp)^m+(pq)^{m-1}= (qp)^{m-1}q+(pq)^{m-1}p$, then $\mathrm{alg}(p,q)\cong  Z_{4m-5}\oplus W_3$ and $\mathrm{dim}\text{ }\mathcal{N}\left(\mathrm{alg}(p,q)\right)=4m-5.$
        \item \label{1k}If $(qp)^{m-1}= (qp)^{m}$ and $(qp)^{m-1}+(pq)^{m-1}\neq (qp)^{m-1}q+(pq)^{m-1}p$ then here $\mathrm{alg}(p,q)\cong  Z_{4m-6}\oplus W_4$ and $\mathrm{dim}\text{ }\mathcal{N}\left(\mathrm{alg}(p,q)\right)=4m-5.$
        \item \label{2k}In the case when $(qp)^{m-1}\neq (qp)^{m}$ and $(qp)^m+(pq)^{m-1}\neq (qp)^{m-1}q+(pq)^{m-1}p$, then here $\mathrm{alg}(p,q)\cong  Z_{4m-5}\oplus W_4$ and $\mathrm{dim}\text{ }\mathcal{N}\left(\mathrm{alg}(p,q)\right)=4m-4.$

    \end{enumerate}
    \begin{proof}
    Since $p\text{ and }q$ satisfies the assumptions of Lemma \ref{inde}; therefore they fit precisely into one of the four scenarios described in Lemma \ref{alg_class}. The complete characterization of the algebra spanned by $p$ and $q$ can be achieved by assessing their defining relations. Therefore, two algebras falling under the same category of Lemma \ref{alg_class} are isomorphic.
        \begin{enumerate}[label=(\roman*)]
            \item Choose $p_0,q_0$ are to be idempotents of the type $Z_{4m-6}$, then $(q_0p_0)^{m-1}=(p_0q_0)^{m-1}=0$. Furthermore let $p,q$ are idempotents of type $W_3.$ Put $$\bf{p}=\begin{bmatrix}
                p_0 & \\
                 & p\\
            \end{bmatrix}\text{ and }\bf{q}=\begin{bmatrix}
                q_0& \\
                 & q\\
            \end{bmatrix},$$ then $\textbf{p},\textbf{q}\in Z_{4m-6}\oplus W_3$ and \begin{align*}
                (\textbf{qp})^{m-1}=\begin{bmatrix}
                    0 & \\
                     & (qp)^{m-1}
                \end{bmatrix}\text{ and }
                (\textbf{pq})^{m-1}=\begin{bmatrix}
                    0 & \\
                     & (pq)^{m-1}
                \end{bmatrix}. 
            \end{align*} Now since $p,q$ are idempotents of type $W_3$, therefore $qpq=q,pqp=p$ and $q+p=qp+pq.$ By using these properties of $p,q$ we get $(\textbf{pq})^m=(\textbf{pq})^{m-1}$, $(\textbf{qp})^m=(\textbf{qp})^{m-1}$ and $$(\textbf{qp})^{m-1}+(\textbf{pq})^{m-1}= (\textbf{pq})^{m-1}\textbf{p}+(\textbf{qp})^{m-1}\textbf{q}.$$ Clearly by the construction of $\textbf{p}$ and $\textbf{q}$ they are not tightly coupled and $m$ is least positive integer satisfying $(\textbf{pq})^{m-1}=(\textbf{pq})^{m}$ and $(\textbf{pq})^{m-2}\textbf{p}\neq (\textbf{pq})^{m-1}\textbf{p}$. Hence in this case alg$(\textbf{p,q})$$\cong  Z_{4m-6}\oplus W_3$ and $$\text{dim }\mathcal{N}\left(\text{alg}(\textbf{p,q})\right)=\text{dim }\mathcal{N}(Z_{4m-6})+\text{dim }\mathcal{N}(W_3)=4m-8+2=4m-6.$$
        \end{enumerate}
        Similarly, like (\ref{kk}), one can prove the other three case (\ref{1k}), (\ref{2k}) and (\ref{3k}).
    \end{proof}
\end{theorem}   
\section{Group and Drazin invertibility}\label{g&d}
One of the main reasons behind the study of the finite-dimensional algebra generated by two tightly coupled idempotents is to understand group and Drazin invertibility in these types of algebras. Since every element in a finite-dimensional algebra is Drazin invertible, group invertibility, in particular properly group invertibility (i.e. when an element is group invertible without being invertible), is the main concern here. The following theorem of \cite{bottcher2012group} addresses this problem of properly group invertibility in $Z_m.$ 
\begin{theorem}\cite{bottcher2012group}\label{wrong2}
    For $m\geq 1$, $I\in Z_m$. Let $A\in Z_m$, such that $$ A=x_1p+y_1q+x_2pq+y_2qp+\cdots,$$ where $\{x_1,x_2\cdots,y_1,y_2,\cdots\}$ are scalars and $p,q$ are idempotents of type $Z_m.$ Then $A$ is properly group invertible if and only if
    \begin{enumerate}[label=(\roman*)]
        \item $A=0;$\\or
        \item either $x_1\neq0,y_1=0$ or $x_1=0,y_1\neq0$ and $0$ is a root of $\psi(t)$ with multiplicity at least $l(m,y_1);$ where \begin{equation*}
    l(m,y_1)=\begin{cases}
        \lceil\frac{m}{4}\rceil-1\text{ if }m=1\text{ }\mathrm{mod}\text{ }4 \text{ and }y_1=0\\
        \lceil\frac{m}{4}\rceil\text{ otherwise }
    \end{cases}, 
    \end{equation*} and $\psi(t)$ is defined in $(\ref{hom})$. Moreover $\sigma(A)=\{x_1,y_1\}$.
    \end{enumerate}
\end{theorem}
Here in Theorem \ref{wrong2}, the claim that the identity element $I\in Z_m$ is not true in general, and due to this, $A$ can be properly group invertible mean while $x_1y_1\neq 0.$ This can be seen in the following example. 
\begin{example}\label{example}
    Let $$p_1=\begin{bmatrix}
        1 & 1 & 0\\
        0 & 0 & 0\\
        0 & 0& 0\\
    \end{bmatrix}\text{ and } q_1=\begin{bmatrix}
        0 & 0 & 0\\
        0& 1&0\\
        0&0&0\\
    \end{bmatrix}\in M_{3}(\mathbb{C}).$$ It is easy to check that $p_1$ and $q_1$ are two idempotents satisfying $q_1p_1=0$. Therefore $p_1,q_1$ are idempotents of the type  $Z_3.$ It is easily seen that rank $(p_1+q_1)=$ rank $(p_1+q_1)^2,$ hence $p_1+q_1$ is properly group invertible. Clearly, here the identity $I\notin\mathrm{alg}(p_1,q_1)$ and the unit element in $\mathrm{alg}(p_1,q_1)$ is the element $$\begin{bmatrix}
        1&0&0\\
        0&1&0\\
        0&0&0\\
    \end{bmatrix}.$$
\end{example}
By \cite[Lemma 4.2]{bottcher2013classification}\begin{equation}
    p+q-pq-qp+pqp+qpq-\cdots
\end{equation} is the unit in $Z_m,$ we will denote it by $I_{Z_m}$.
 Due to the preservation of group invertibility under isomorphism, B{\"o}ttcher and Spitkovsky \cite[Theorem 2.2]{bottcher2012group} considered only the case when algebra $Z_m$ is generated by \begin{equation}
    p=\begin{bmatrix}\label{wrong}
        I&B\\
        0&0\\
    \end{bmatrix} \text{ and }q=\begin{bmatrix}
        0&0\\
        C&I\\
    \end{bmatrix},
\end{equation} for some suitable $B$ and $C$. For the $p,q$ given in (\ref{wrong}), $I=I_{Z_m}$. Now if we consider the idempotents 
\begin{equation}
    p^{\prime}=\begin{bmatrix}\label{right}
        I&B_1&B_2\\
        0&0&0\\
        0&0&0\\
    \end{bmatrix} \text{ and }q^{\prime}=\begin{bmatrix}
        0&0&0\\
        C_1&I&C_2\\
        0&0&0\\
    \end{bmatrix},
\end{equation} then for suitable $B_1,C_1$, alg$(p^{\prime},q^{\prime})=Z_m$, but here $I_{Z_m}\neq I$. 
Note that similarly as given in \cite{bottcher2012group}, depending on $m$, it is possible to choose $B_1,C_1$ such that $p^{\prime},q^{\prime}$ are of type $Z_m$, and no additional condition is required on $B_2,C_2$.
\\ Now if $Z_m=\text{alg}(p,q)$ then any elements $A\in Z_m$ has the representation \begin{equation}\label{represe}
    A=x_1p+y_1q+x_2pq+y_2qp+\cdots,
\end{equation} where $x_i,y_i$ are scalars. For further use, we recall the following functions introduced by \cite{bottcher2012group}
\begin{align}
    \upvarphi_{00}(t)&=x_1+(x_2+x_3)t+(x_4+x_5)t^2+\cdots\notag\\
    \upvarphi_{11}(t)&=y_1+(y_2+y_3)t+(y_4+y_5)t^2+\cdots\notag\\
    \upvarphi_{01}(t)&=(x_1+x_2)+(x_3+x_4)t+\cdots \label{hom}\\
    \upvarphi_{10}(t)&=(y_1+y_2)+(y_3+y_4)t+\cdots\notag\\
    \psi(t)&=\upvarphi_{00}(t)\upvarphi_{11}(t)-t\upvarphi_{01}(t)\upvarphi_{10}(t)\notag,
\end{align} and we also define \begin{align}
    \upvarphi_{02}(t)&=x_1+x_3t+x_5t^2+\cdots\notag\\
    \upvarphi_{12}(t)&=y_1+y_3t+y_5t^2+\cdots\notag\\
    \upvarphi_{02^{\prime}}(t)&=x_2+x_4t+x_6t^2+\cdots\label{hom2}\\
    \upvarphi_{12^{\prime}}(t)&=y_2+y_4t+y_6t^2+\cdots\notag\\
    \psi_1(t)&=\upvarphi_{00}(t)\upvarphi_{12}(t)-t\upvarphi_{10}(t)\upvarphi_{02^{\prime}(t)}\notag\\ 
    \psi_{2}(t)&=\upvarphi_{00}(t)\upvarphi_{12^{\prime}}(t)-\upvarphi_{10}(t)\upvarphi_{02}(t)\notag
\end{align}
Our next lemma establishes a relation between zeros of the polynomials $\psi,\psi_1$ and $\psi_2$ when $y_1=0$.
\begin{lemma}\label{countzero}
Whenever $y_1=0,$ in $(\ref{hom})$ and $(\ref{hom2})$,
    if $0$ is a root of the polynomial $\psi$ with multiplicity $n$, then $0$ is also a root of the polynomials $\psi_1$ and $\psi_2$ with at least the same multiplicity.
    \begin{proof}
        Since $y_1=0$ then we obtain $$\psi(0)=\psi_1(0)=\psi_2(0)=0.$$ Thus, the statement holds for $n=1$ case. Let assume it is true for $n=n-1$ case, i.e. if $0$ is a root of the polynomial $\psi$ with multiplicity $n-1$, then $0$ is also a root of $\psi_1$ and $\psi_2$ with multiplicity $n-1.$ Now assume $0$ as a root of $\psi$ with multiplicity $n$, then $$\psi(0)=\diff*{\psi}{t}{(0)}=\cdots=\diff*[n]{\psi}{t}{(0)}=0.$$ Again from $\diff*[n]{\psi}{t}{(0)}=0$, choosing $x_k,y_k=0$ for $k< 1$ we obtain
        \begin{align*}
            \diff*[n]{\psi}{t}{(0)}=&\sum_{k=0}^{n}\left((x_{2n-2k}+x_{2n-2k+1})(y_{2k}+y_{2k+1})\right.\\{}&-\left.(x_{2n-2k-1}+x_{2n-2k})(y_{2k+1})(y_{2k+2})\right)\\=&-x_{2n}y_2+x_{2n-1}y_3+x_{2n-2}(y_2-y_4)+x_{2n-3}(y_5-b_3)\\{}&+\cdots +x_2(y_{2n-2}-y_{2n})+x_1(y_{2n+1}-y_{2n-1})\\=&{}0.
        \end{align*}
        
        By previous assumption $0$ is also a root of $\psi_1$ with multiplicity $n-1$, therefore  \begin{align*}
            \diff*[n-1]{\psi_1}{t}{(0)}=&{}\sum_{k=0}^{n-1}\left((x_{2n-2k-2}+x_{2n-2k-1})y_{2k+1}\right.\\{}&-\left.(y_{2n-2k-3}+y_{2n-2k-2})x_{2k+2}\right)\\=&{}-x_{2n-2}y_2+x_{2n-3}y_3-x_{2n-4}y_4+\cdots -x_2y_{2n-2}+x_1y_{2n-1}\\{}=&0.
        \end{align*} 
        Now \begin{align*}
            \diff*[n]{\psi_1}{t}{(0)}=&{}\sum_{k=0}^{n}\left((x_{2n-2k}+x_{2n-2k+1})y_{2k+1}\right.\\{}&-\left.(y_{2n-2k-1}+y_{2n-2k})x_{2k+2}\right)\\=&-x_{2n}y_2+x_{2n-1}y_3-x_{2n-2}y_4+\cdots -x_2y_{2n}+x_1y_{2n+1}\\=&
            \diff*[n]{\psi}{t}{(0)}+\diff*[n-1]{\psi_1}{t}{(0)}.
        \end{align*}
       Therefore by $\diff*[n]{\psi}{t}{(0)}=0\text{ and }\diff*[n-1]{\psi_1}{t}{(0)}=0$ it follows that $ \diff*[n]{\psi_1}{t}{(0)}=0.$ Hence $0$ is a root of the function $\diff*[n]{\psi_1}{t}{}$, it implies $0$ is a root of $\psi_1$ with multiplicity at least $n.$ Similarly one can proof for $\psi_2$ also.
    \end{proof}
\end{lemma}
Now we are ready to present the version of \cite[Theorem 2.2]{bottcher2012group}, when $I\notin Z_m.$
\begin{theorem}\label{Z_m}
    Let $n\geq 1$, $I\notin Z_m$ and $p,q$ are idempotent of type $Z_m.$ If $A\in Z_m$ is in the form $(\ref{represe})$, then $A$ is properly group invertible if and only if one of the following is true
    \begin{enumerate}[label=(\roman*)]
        \item \label{11k}$A=0;$
         \item \label{13k} $x_1y_1\neq0$;
        \item \label{12k} either $x_1\neq 0,y_1=0\text{ or }x_1=0,y_1\neq0$ and $0$ is a root of $\psi(t)$ with multiplicity at least $\lceil\frac{m}{4}\rceil.$ 
    \end{enumerate} Moreover $\sigma(A)=\{x_1,y_1,0\}.$ 
    \begin{proof}
      The case when $m=1$ or $2$ is trivial, so let $m\geq 3.$  It is straightforward to verify that if $p,q$ is of the form (\ref{right}), then $I\notin Z_m$ and
        \begin{equation*}
           A=\begin{bmatrix}
               \upvarphi_{00}(B_1C_1)&\upvarphi_{01}(B_1C_1)B_1&\upvarphi_{02}(B_1C_1)B_2+\upvarphi_{02^{\prime}}(B_1C_1)B_1C_2\\
               C_1\upvarphi_{10}(B_1C_1)&\upvarphi_{11}(C_1B_1)&\upvarphi_{12}(C_1B_1)C_2+\upvarphi_{12^{\prime}}(C_1B_1)C_1B_2\\
               0&0&0\\
           \end{bmatrix}. 
        \end{equation*}
        As mentioned in \cite{bottcher2012group} here the matrix $B_1C_1$ is nilpotent with degree $\lceil\frac{m}{4}\rceil-1$ when $m=1\text{ mod }4$ and $B_1C_1,\text{ }C_1B_1$ are nilpotent with degree $\lceil\frac{m}{4}\rceil$ in other cases.
        Now by \cite[Theorem 3.1]{meyer1977index}, if the $2\times 2$ block\begin{equation}\label{111}
           A_1= \begin{bmatrix}
                 \upvarphi_{00}(B_1C_1)&\upvarphi_{01}(B_1C_1)B_1\\
               C_1\upvarphi_{10}(B_1C_1)&\upvarphi_{11}(C_1B_1)\\
            \end{bmatrix}
        \end{equation} is invertible, then $A$ becomes properly group invertible. According to Theorem \ref{wrong2}, if $x_1y_1\neq 0$, then $A_1$ is invertible; hence for the case $x_1y_1\neq 0,$ $A$ is properly group invertible. Next, if $x_1=y_1=0$, then $A$ is a nilpotent matrix, which is group invertible if and only if $A=0.$ Thus $(\ref{11k})$ and $(\ref{13k})$ follows. 
       \\ Now assume $x_1\neq 0, y_1=0.$ Since $\upvarphi_{00}(0)=x_1\neq 0$ and $B_1C_1$ is nilpotent, therefore $\upvarphi_{00}(B_1C_1)$ is invertible, and the Schur complement of $\upvarphi_{00}(B_1C_1)$ in $A$ is 
        \begin{equation}\label{schur}
           \frac{1}{\upvarphi_{00}(t)}\begin{bmatrix}
                {\psi(t)}(C_1B_1)&{\psi_1(t)}(C_1B_1)C_2+{\psi_2(t)}(C_1B_1)C_1B_2\\
                0&0\\
            \end{bmatrix}.
        \end{equation} But $\psi(0)=0$, so (\ref{schur}) is a nilpotent matrix. Again the Schur complement of $\upvarphi_{00}^2(t)-t\upvarphi_{01}(t)\upvarphi_{10}(t)(B_1C_1)$ in $A^2$ is  
        \begin{equation}\label{schur2}\notag
           \frac{1}{\upvarphi_{00}^2(t)-t\upvarphi_{01}(t)\upvarphi_{10}(t)}\begin{bmatrix}
                {\psi(t)}(C_1B_1)&{\psi_1(t)}(C_1B_1)C_2+{\psi_2(t)}(C_1B_1)C_1B_2\\
                0&0\\
            \end{bmatrix}^2. 
        \end{equation} Thus rank $(A)$=rank $(A^2)$ if and only if Schur complement of $\upvarphi_{00}(t)(B_1C_1)$ in $A$ and Schur complement of $\upvarphi_{00}^2(t)-t\upvarphi_{01}(t)\upvarphi_{10}(t)(B_1C_1)$ in $A^2$ have the same rank. Since $(\ref{schur})$ is a nilpotent matrix; therefore, this is possible only when $(\ref{schur})$ is the zero matrix. Now for $Z_m$, $B_1C_1$ is nilpotent matrix of degree $\lceil\frac{m}{4}\rceil-1$ when $m=1$ mod $4$ and $B_1C_1,C_1B_1$ are nipotent with degree $\lceil\frac{m}{4}\rceil$ otherwise. Hence ${\psi(t)}(C_1B_1)=0$ if and only if  $0$ is a root of $\psi(t)$ with multiplicity atleast $\lceil\frac{m}{4}\rceil$. According to Lemma \ref{countzero}, if $0$ is a root of $\psi(t)$ with multiplicity $\lceil\frac{m}{4}\rceil$, then $0$ is also a root of $\psi_1(t)$ and $\psi_2(t)$ with at least the same multiplicity.  Hence we get ${\psi_1(t)}(C_1B_1)C_2+{\psi_2(t)}(C_1B_1)C_1B_2=0$. This proves $(\ref{12k})$ when $x_1\neq 0$ and $y_1=0.$ Likewise, the situation  when $x_1=0$ but $y_1\neq 0$ can be treated by considering the Schur complement of $\upvarphi_{11}(C_1B_1).$ Moreover since $\sigma(A)=\sigma(A_1)\cup \{0\},$ therefore from Theorem \ref{wrong2} it also follows that $\sigma(A)=\{0,x_1,y_1\}$.
    \end{proof}   
\end{theorem}
 With the help of the previous theorem and \cite[Theorem 3.5]{bottcher2012group}, it is possible to classify all the properly group invertible elements in the algebra $Z_m\oplus W_3$ and $Z_m\oplus W_4.$ Our next two results describe that.
\begin{theorem}
    If $A\in Z_m\oplus W_3$ is of the form $(\ref{represe}),$ then $A$ is properly group invertible if and only if one of the following is true
    \begin{enumerate}[label=(\roman*)]
        \item $A=0$;
        \item $\displaystyle\sum_i(x_{2i-1}+y_{2i})=\sum_i(x_{2i}-y_{2i})=\sum_iy_{i}=0$ and \begin{enumerate}
            \item $x_1y_1\neq 0$;  \\ or \item  $x_1,y_1$ satisfies condition (\ref{12k}) of Theorem $\ref{Z_m}$;
        \end{enumerate}
        \item $\sum_i (x_i+y_i)\neq 0$ and 
        \begin{enumerate}
            \item all the coefficients of $A$ in $(\ref{represe}) $ are $0$, except for the last three;\\ or 
            \item $x_1y_1\neq 0$; \\or
             \item $x_1,y_1$ satisfies condition (\ref{12k}) of Theorem $\ref{Z_m}$.
        \end{enumerate}
    \end{enumerate} Moreover $\sigma(A)=\{0,x_1,y_1,x+y\} $ where $x=\sum_{i}x_i$ and $y=\sum_iy_i.$
    \begin{proof}
        Let $\textbf{p},\textbf{q}$ are idempotents of type $Z_m\oplus W_3$  then\begin{equation} \label{formation}
            \bf{p}=\begin{bmatrix}
                p_0&\\ &p
            \end{bmatrix}\text{ and } \bf{q}=\begin{bmatrix}
                q_0&\\&q
            \end{bmatrix}
        \end{equation} where $p_0,q_0$ are type $Z_m$ idempotents and $p,q$ are idempotents of type $W_3.$ Denote $A_0=x_1p_0+y_1q_0+x_2p_0q_0\cdots$ and $A_1=x_1p+y_1q+x_2pq\cdots.$ Then by Theorem \ref{Z_m} and \cite[Theorem 3.5]{bottcher2012group} $$\sigma(A)=\sigma(A_0)\cup \sigma(A_1)=\{0,x_1,y_1,x+y\}.$$ According to \cite[Theorem 3.5]{bottcher2012group} $A_1$ is not invertible; therefore $A$ is properly group invertible if and only if $A_0$ and $A_1$ are group invertible. Since in $W_3$ 
        \begin{align*}
            A_1&=x_1p+y_1q+x_2pq\cdots\\
            &=\sum_i(x_{2i-1}+y_{2i})p+\sum_i(x_{2i}-y_{2i})pq+\sum_i(y_{2i-1}+y_{2i})q,
        \end{align*} hence by \cite[Theorem 3.5]{bottcher2012group} $A_1$ is group invertible if and only if one of the following holds\begin{enumerate}
            \item $x_i,y_i=0$ for all $i$;
            \item $\displaystyle\sum_i(x_{2i-1}+y_{2i})=\sum_i(x_{2i}-y_{2i})=\sum_i(y_{2i-1}+y_{2i})=0$;
            \item $\displaystyle\sum_i x_i+y_i\neq 0$.
        \end{enumerate} From here, the desired result can be confirmed with the help of Theorem \ref{Z_m}.
    \end{proof}
\end{theorem}
\begin{theorem}
    If $A\in Z_m\oplus W_4$ is of the form $(\ref{represe})$ then $A$ is properly group invertible if one of the followings holds
    \begin{enumerate}[label=(\roman*)]
        \item $A=0$;
        \item $\displaystyle\sum_ix_{2i-1}=\sum_ix_{2i}=\sum_iy_{2i-1}=\sum_iy_{2i}=0$ and \begin{enumerate}
            \item $x_1y_1\neq 0$; \\ or \item  $x_1,y_1$ satisfies condition (\ref{12k}) of Theorem $\ref{Z_m}$;;
        \end{enumerate}
        \item $\displaystyle\sum_i( x_i+y_i)\neq 0$, $\displaystyle\left(\sum_ix_{2i-1}\right)\left(\sum_iy_{2i-1}\right)=\left(\sum_ix_{2i}\right)\left(\sum_iy_{2i}\right)$ and 
        \begin{enumerate}
            \item all the coefficients of $A$ in $(\ref{represe})$ are $0$, except for the last four;\\ or 
            \item $x_1y_1\neq 0$;\\ or
             \item  $x_1,y_1$ satisfies condition (\ref{12k}) of Theorem $\ref{Z_m}$.
        \end{enumerate}
    \end{enumerate} Moreover $\sigma(A)=\{0,x_1,y_1,x+y\} $ where $x=\sum_{i}x_i$ and $y=\sum_iy_i.$
    \begin{proof}
        This follows from Theorem \ref{Z_m} and \cite[Theorem 3.5]{bottcher2012group}.
    \end{proof}
\end{theorem}
As already mentioned, being a finite-dimensional algebra, every element in $Z_m$ is Drazin invertible, and it is possible to provide an upper bound for the Drazin index of any element $A\in Z_m.$ The following theorem describes this.
\begin{theorem}\label{index}
    Let $m\geq 1$, and $A\in Z_m$ is not group invertible. If $A$ is nilpotent then $$i(A)\leq\begin{cases}
       2 \lfloor\frac{m}{4}\rfloor\text{ if }m=0 \text{ }\mathrm{mod}\text{ }4\\
       2 \lfloor\frac{m}{4}\rfloor+1\text{ if }m=1\text{ or }2 \text{ }\mathrm{mod}\text{ }4\\
        2 \lfloor\frac{m}{4}\rfloor+2\text{ if }m=3 \text{ }\mathrm{mod}\text{ }4
    \end{cases},$$ and if $A$ is not nilpotent then $i(A)\leq \lceil\frac{m}{4}\rceil.$
    \begin{proof}
        Since $A\in Z_m$, therefore $A$ is of the form (\ref{represe}). Now if $x_1y_1\neq0,$ then by Theorem \ref{Z_m}, $A$ is group invertible, and if exactly one of $x_1,y_1$ is $0$, then by Theorem \ref{Z_m}, $A$ is Drazin invertible with Drazin index $\leq \lceil\frac{m}{4}\rceil$. Finally if $x_1=0=y_1$, then $$A=x_2pq+y_2qp+x_3pqp+y_3qpq+\cdots,$$ which is a nilpotent element in alg$(p,q)$. Here using induction, it becomes evident that $ (pq)_{k+1}$ and $(qp)_{k+1}$ are the elements with the lowest order in the representation of $A^k.$ If $m=0$ mod $4$, then any element in $Z_m$ with order $\geq  2 \lfloor\frac{m}{4}\rfloor+1$ is equals to $0,$ hence $A^{ 2 \lfloor\frac{m}{4}\rfloor}=0.$ Therefore $i(A)\leq  2 \lfloor\frac{m}{4}\rfloor.$ Likewise, one can confirm the validity of other cases as well.
    \end{proof}
\end{theorem}
\begin{remark}
    Using Theorem $\ref{index}$, one can find an upper for the Drazin index of elements in algebra $U_m,D_m,D_m^*$.
\end{remark}
\section{Representation of Drazin and group inverse}\label{representation}
Given that $p,q$ satisfies $(pq)^{m-1}=(pq)^{m}$, then the finite dimensionality of the algebra alg$(p,q)$ ensures that the Drazin inverse of $(\alpha p+q)$  always exists and lies in alg$(p,q)$.  However, expressing the representation of $(\alpha p+q)^D$ using $p,q$ is still a challenging task. The representation of $(\alpha p+\beta q)$ under the assumption $(pq)^n=(pq)^{n-1}$ is already provided by Shi and Guolin \cite{shi2013drazin}  for idempotents in a Banach algebra, but they proved it using the integral representation of the Drazin inverse in a Banach algebra. But here, we provide a proof which is also suitable for an associative algebra.

\begin{theorem}\label{p,q}
    Let $\alpha$ be a non-zero scalar. If $(pq)^{m-1}=(pq)^{m}$ then $(\alpha p+q)$ is Drazin invertible and $i(\alpha p+q)\leq \begin{cases}
        2, \text{ if }\alpha\neq -1\\
        3, \text{ if }\alpha=-1
    \end{cases}.$  Moreover 
    \begin{equation}\label{fc}
        (\alpha p+q)^D=\begin{cases}
        \textbf{A}^3(\alpha p+q)^2,\text{ if }\alpha\neq-1\\
        \textbf{B}^4(\alpha p+q)^3,\text{ if }\alpha=-1
    \end{cases},
    \end{equation}
     where 
    \begin{align*}
     \textbf{A}=&\displaystyle\sum_{i=1}^{2m-3}(-1)^{i-1}\left(\left(\lfloor \frac{i}{2}\rfloor+\frac{\lfloor\frac{i}{2}\rfloor+\phi(i)}{\alpha}\right)(pq)^{\lfloor \frac{i}{2}\rfloor}p^{\phi(i)}\right.\\{}&+\left.\left(\lfloor \frac{i}{2}\rfloor+\phi(i)+\frac{\lfloor\frac{i}{2}\rfloor}{\alpha}\right)(qp)^{\lfloor \frac{i}{2}\rfloor}q^{\phi(i)}\right)\\{}&-\left(m-1+\frac{m-1}{\alpha}\right)(qp)^{m-1}-\left(\frac{m-1}{\alpha}-\frac{\alpha}{(1+\alpha)^2}\right)(pq)^{m-1}\\{}&+\left(\frac{m-1}{\alpha}+\frac{1}{(1+\alpha)^2}\right)(qp)^{m-1}q,
    \end{align*} and $$\textbf{B}=\displaystyle\sum_{i=1}^{2m-3}\phi(i)\left((qp)^{\lfloor \frac{i}{2}\rfloor}q^{\phi(i)}-(pq)^{\lfloor \frac{i}{2}\rfloor}p^{\phi(i)}\right)$$ here the function $\phi$ is defined as $\phi(i)=\begin{cases}
        0, \text{ if }i\text{ is even,}\\1, \text{ if }i
\text{ is odd}    
      \end{cases}$.
      \begin{proof}
           First let assume that $\alpha\neq -1$, then by \cite[Lemma]{puystjens2004drazin} to prove (\ref{fc}), it is enough to prove that $
           \textbf{A}(\alpha p+q)^3=(\alpha p+q)^2$ and $(\alpha p+q)^3\textbf{A}^{\prime}=(\alpha p+q
           )^2,$ for some $\textbf{A}^{\prime}\in\mathcal{A}.$ Now we compute the coefficient of $(pq)^{\lfloor\frac{i}{2}\rfloor}p^{\phi(i)}$ in the expansion $\textbf{A}(\alpha p+q)^3$ for $2<i\leq 2m-3$. Let assume $i$ to be even, then $i=2k$ for some $k\in\mathbb{N}$, and the coefficient of $(pq)^{\lfloor\frac{i}{2}\rfloor}p^{\phi(i)}\left(=(pq)^k\right)$ is \begin{align*}
               \left(\frac{k-1}{\alpha}+k-2\right)&\alpha-\left(\frac{k-1}{\alpha}+k-1\right)(2\alpha+\alpha^2)\\&+\left(\frac{k}{\alpha}+k-1\right)(1+\alpha+\alpha^2)-\left(\frac{k}{\alpha}+k\right)\\{}=&0.
           \end{align*} Next if $i$ is odd then $i=2k+1$ for some $k\in\mathbb{N},$ then the coefficient of $(pq)^{\lfloor\frac{i}{2}\rfloor}p^{\phi(i)}\left(=(pq)^kp\right)$ is \begin{align*}
               -\left(\frac{k-1}{\alpha}+k-1\right)&\alpha^2+\left(\frac{k}{\alpha}+k-1\right)(2\alpha^2+\alpha)\\{}&-\left(\frac{k}{\alpha}+k\right)(\alpha+\alpha^2+\alpha^3)+\left(\frac{k+1}{\alpha}+k\right)\alpha^3\\{}=& 0. 
           \end{align*} In particular, the coefficient of $p\text{ and }pq$ in $\textbf{A}(\alpha p+q)^3$ are $\alpha^2$ and $\alpha$, respectively. Hence, we obtain that in the expansion of $\textbf{A}(\alpha p+q)^3$, the coefficient of $(pq)^{\lfloor \frac{i}{2}\rfloor}p^{\phi(i)}$ is $$\begin{cases}
            \alpha^2, \text{ if }i=1\\
            \alpha, \text{ if }i=2\\
            0, \text{ if }3\leq i\leq 2m-3
        \end{cases}.$$ Likewise, it can be verified that the coefficient of $(qp)^{\lfloor \frac{i}{2}\rfloor}q^{\phi(i)}$ is $$\begin{cases}
            1, \text{ if }i=1\\
            \alpha, \text{ if }i=2\\
            0, \text{ if }3\leq i\leq 2m-2
        \end{cases}.$$Now, if $i=2m-2$, then the coefficient of $(pq)^{\lfloor \frac{i}{2}\rfloor}p^{\phi(i)}\left(=(pq)^{m-1}\right)$ is 
        \begin{align*}
            \left(\frac{m-2}{\alpha}+m-3\right)&\alpha-\left(\frac{m-2}{\alpha}+m-2\right)(2\alpha+\alpha^2)\\{}&+\left(\frac{m-1}{\alpha}+m-2\right)(1+2\alpha+\alpha^2)\\{}&-\left(\frac{m-1}{\alpha}-\frac{\alpha}{(1+\alpha)^2}\right)(1+\alpha)^2\\=&0
        \end{align*}In a similar manner, it can be confirmed that the coefficient of $(pa)^{m-1}p,$ $(qp)^{m-1}q$ and $(qp)^m$ is $0$ in the expansion of $\textbf{A}(\alpha p+q)^3.$ Hence, through evaluating the coefficients of each term in the expansion of $\textbf{A}(\alpha p+q)^3$, we obtain that $$\textbf{A}(\alpha p+q)^3=(\alpha p+q)^2.$$ Certainly, by choosing \begin{align*}
     \textbf{A}^{\prime}=&\displaystyle\sum_{i=1}^{2m-3}(-1)^{i-1}\left(\left(\lfloor \frac{i}{2}\rfloor+\frac{\lfloor\frac{i}{2}\rfloor+\phi(i)}{\alpha}\right)(pq)^{\lfloor \frac{i}{2}\rfloor}p^{\phi(i)}\right.\\{}&\left.+\left(\lfloor \frac{i}{2}\rfloor+\phi(i)+\frac{\lfloor\frac{i}{2}\rfloor}{\alpha}\right)(qp)^{\lfloor \frac{i}{2}\rfloor}q^{\phi(i)}\right)+\left(m-1+\frac{m-1}{\alpha}\right)(qp)^{m-1}\\{}&+\left(\frac{1}{(\alpha+1)^2}+1-m\right)(pq)^{m-1}+\left(m-1+\frac{\alpha}{(1+\alpha)^2}\right)(pq)^{m-1}p,
    \end{align*} one can proof $(\alpha p+q)^3\textbf{A}^{\prime}=(\alpha p+q)^2.$ Therefore if $\alpha\neq -1$ then $i(\alpha p+q)\leq 2$ and $$(\alpha p+q)^D=\textbf{A}^3(\alpha p+q)^2.$$\\
    Now if $\alpha=-1$ then\begin{align*}
        \textbf{B}(q-p)^4=&\left(\displaystyle\sum_{i=1}^{2m-3}\phi(i)\left((qp)^{\lfloor \frac{i}{2}\rfloor}q^{\phi(i)}-(pq)^{\lfloor \frac{i}{2}\rfloor}p^{\phi(i)}\right)\right)(q-p)^4\\{}=&\left((q-p)+(pq)^{m-1}p-(qp)^{m-1}q\right)(q-p)^2\\{}=&(q-p)^3+\left((qp)^{m}q-(qp)^{m-1}q\right.\\{}&\left.+(pq)^{m-1}q-(pq)^mp\right)\\{}=&{(q-p)^3}\\{}=&(q-p)^4\textbf{B}.
    \end{align*} Therefore $i(q-p)\leq 3$ and $(q-p)^D=\textbf{B}^4(\alpha p+q)^3.$
      \end{proof}
\end{theorem}
In this situation, the next question is if $\lambda(pq)^{m-1}=(pq)^m$ for some non-unit scalar $\lambda$, then what conclusions can be drawn regarding the Drazin invertibility of $(\alpha p+q).$ Chen et al. \cite{chen2020drazin} recently investigated a specific instance of this problem, specifically when $pqp=\lambda p.$ In their study, they  established the group invertibility $c_1p+c_2q+c_3pq+c_4qp+c_5pqp,$ where $c_1,\cdots,c_5$ are scalar. In the subsequent theorem, we establish the group invertibility of $(\alpha p+q),$ when $(pq)^m=\lambda(pq)^{m-1}.$
\begin{theorem}
    Let $\alpha$ be a non-zero scalar. If $\lambda(pq)^{m-1}=(pq)^m$ for some  non-unit scalar $\lambda,$ then $(\alpha p+q)$ is group invertible and \begin{equation}\label{case}
        (\alpha p+q)^g=\begin{cases}
        \textbf{A}^2(\alpha p+q),\text{ if }\alpha\neq-1\\
        \textbf{B}^2(\alpha p+q),\text{ if }\alpha=-1
    \end{cases},
    \end{equation} where 
    \begin{align*}
     \textbf{A}=&\displaystyle\sum_{i=1}^{2m-3}(-1)^{i-1}\left(\left(\lfloor \frac{i}{2}\rfloor+\frac{\lfloor\frac{i}{2}\rfloor+\phi(i)}{\alpha}\right)(pq)^{\lfloor \frac{i}{2}\rfloor}p^{\phi(i)}\right.\\{}&\left.+\left(\lfloor \frac{i}{2}\rfloor+\phi(i)+\frac{\lfloor\frac{i}{2}\rfloor}{\alpha}\right)(qp)^{\lfloor \frac{i}{2}\rfloor}q^{\phi(i)}\right)-\left(m-1+\frac{m-1}{\alpha}\right)(qp)^{m-1}\\{}&+a_1(pq)^{m-1}+a_2(pq)^{m-1}p+b_1(qp)^{m-1}q+b_2(qp)^m,
    \end{align*} and $$\textbf{B}=\displaystyle\sum_{i=1}^{2m-3}\phi(i)\left((qp)^{\lfloor \frac{i}{2}\rfloor}q^{\phi(i)}-(pq)^{\lfloor \frac{i}{2}\rfloor}p^{\phi(i)}\right)+\frac{(qp)^{m-1}q-(pq)^{m-1}p}{1-\lambda}$$ here $\phi$ is the same function as defined in Theorem $\ref{p,q}$, and\begin{align*}
        a_1={}&\frac{(\alpha+1)\left(m(\lambda-1)+1-2\lambda\right)}{\alpha(\lambda-1)^2},\\a_2={}&\frac{m(1+\alpha-\lambda-\alpha\lambda)+(\lambda-\alpha+2\alpha\lambda)}{\alpha(\lambda-1)^2},\\b_1={}&\frac{m(\alpha+1)(1-\lambda)+\alpha\lambda+2\lambda-1}{\alpha(\lambda-1)^2},\\b_2={}&-\frac{m(1+\alpha-\alpha\lambda-\lambda)+\lambda(\alpha+1)}{\alpha(\lambda-1)^2}.
    \end{align*}
    \begin{proof}
        First, consider the case when $\alpha\neq -1$. Similar to the proof of Theorem \ref{p,q}, to prove (\ref{case}), it is enough if we prove that $\textbf{A}(\alpha p+q)^2=\alpha p+q$ and $(\alpha p+q)^2\textbf{A}=\alpha p+q.$ Now, one can verify that in the expansion of $\textbf{A}(\alpha p+q)^2$, the coefficient of $(pq)^{\lfloor \frac{i}{2}\rfloor}p^{\phi(i)}$ is $$\begin{cases}
            \alpha, \text{ if }i=1\\
            0, \text{ if }2\leq i\leq 2m-3
        \end{cases}$$ and the coefficient of $(qp)^{\lfloor \frac{i}{2}\rfloor}q^{\phi(i)}$ is $$\begin{cases}
            1, \text{ if }i=1\\
            0, \text{ if }2\leq i\leq 2m-2
        \end{cases}.$$ Next the coefficient of $(pq)^{m-1}$ in the expansion of $\textbf{A}(\alpha p+q)^2$ is \begin{align}\label{coe}
            \notag-\left(\frac{m-2}{\alpha}+m-2\right)\alpha+&\left(\frac{m-1}{\alpha}+m-2\right)(1+\alpha)+a_1(1+\lambda\alpha)\\\notag{}&+a_2(\lambda+\lambda\alpha)\\{}=&(m-1)\left(1+\frac{1}{\alpha}\right)+a_1(1+\alpha\lambda)+a_2(\lambda+\lambda\alpha).
        \end{align}
        Identically we obtain the coefficient of $(pq)^{m-1}p$ in $\textbf{A}(\alpha p+q)^2$ is \begin{align}\label{PQmp}
            \left(\frac{m-1}{\alpha}+m-2\right)\alpha+a_1(\alpha^2+\alpha)+a_2(\alpha^2+\lambda\alpha).
        \end{align}
        Now solving the system of linear equations obtained from (\ref{coe}) and (\ref{PQmp}) i.e. 
        \begin{align*}
            (m-1)\left(1+\frac{1}{\alpha}\right)+a_1(1+\alpha\lambda)+a_2(\lambda+\lambda\alpha)=&0\text{ and }\\
            \left(\frac{m-1}{\alpha}+m-2\right)+a_1(\alpha+1)+a_2(\alpha+\lambda)=&0
        \end{align*} we get 
        \begin{align*}
             a_1={}&\frac{(\alpha+1)\left(m(\lambda-1)+1-2\lambda\right)}{\alpha(\lambda-1)^2},\\a_2={}&\frac{m(1+\alpha-\lambda-\alpha\lambda)+(\lambda-\alpha+2\alpha\lambda)}{\alpha(\lambda-1)^2}.
        \end{align*}
        Analogously solving the system of linear equations obtained from the coefficient of $(qp)^{m-1}q$ and $(qp)^{m}$ in  $\textbf{A}(\alpha p+q)^2$, the values of  $b_1,\text{ }b_2$ can be confirmed. Hence we get $$\textbf{A}(\alpha p+q)^2=\alpha p+q.$$ In the same way, by comparing the coefficient of each term in the expansion of $(\alpha p+q)^2\textbf{A}$ one can verify $$(\alpha p+q)^2\textbf{A}=\alpha p+q.$$ Hence $(\alpha p+q)$ is group invertible and $$(\alpha p+q)^g=\textbf{A}^2(\alpha p+q).$$ \\
        Now if $\alpha=-1$ then \begin{align*}
            \textbf{B}(q-p)^2=&(q-p)+(pq)^{m-1}p-(qp)^{m-1}q+\frac{(qp)^{m-1}q-(pq)^{m-1}p}{1-\lambda}(q-p)^2\\
            =&(q-p)\\=&(q-p)^2\textbf{B}.
        \end{align*} Therefore $(q-p)$ is group invertible and $$(q-p)^g=\textbf{B}^2(q-p).$$
    \end{proof}
\end{theorem}






\end{document}